\documentclass{article}
\usepackage{graphicx} % Required for inserting images

\usepackage{xcolor}
\usepackage{cite}
\usepackage{amsmath,amssymb,amsfonts}
\usepackage{algorithmic}
\usepackage{mathtools}
\usepackage{graphicx}
\usepackage{amsthm}

\graphicspath{{images/}}

\usepackage{textcomp}
\usepackage{xcolor}
\def\BibTeX{{\rm B\kern-.05em{\sc i\kern-.025em b}\kern-.08em
    T\kern-.1667em\lower.7ex\hbox{E}\kern-.125emX}}
\title{Passive Stability and Adaptive Control of Teleoperated System using Wave Variables and Predictor Techniques}
\author{Sridhar B Mudhangulla, Naveen K Rajarajan, Olugbenga Moses Anubi}
\date{}

\begin{document}

\maketitle
\begin{centering}
    Department of Electrical and Computer Engineering, Florida State University
E-mail: nr21w@fsu.edu, sm19ch@fsu.edu, oanubi@fsu.edu
\end{centering}

\section*{Abstract}
This paper addresses the challenge of achieving stable adaptive teleoperation and improving the convergence rate in the presence of high communication time delays. We employ a passivity-based formalism to establish stability using wave variables and wave scattering techniques, and we enhance the convergence rate by combining it with predictor-based approaches. The elevated time delay within the teleoperated communication layer is known to induce an oscillatory behavior, which reduces the convergence rate and increases the settling time in the convergence of power variables. This issue is addressed in this paper by utilizing a Smith predictor on the operator end and Minimum Jerk (MJ) predictor on the remote end. We present experimental and simulation results to demonstrate the improvements, ensuring stable teleoperation under high communication time delays.

\section{Introduction}
In an era marked by advancements in robotics and automation, teleoperation stands as a crucial bridge between human control and autonomous machines \cite{prakash2022vehicle}. In recent years, teleoperation has seen rapid advancements, driven by the increasing demand for remotely controlled robotic systems in various domains such as healthcare, manufacturing, disaster response, and space exploration \cite{telerobotics}. Teleoperation systems offer the promise of extending human capabilities by allowing operators to remotely control robots in environments that are dangerous, inaccessible, or simply distant \cite{Hokayem}. However, despite these remarkable developments, the achievement of stable and intuitive teleoperation remains a formidable challenge, particularly when communication delays are introduced into the control loop \cite{anderson_bilateral}. 

Communication delays are an inherent and inevitable aspect of remote teleoperation systems, arising from the finite speed of radio-frequency (RF) signals as they traverse vast distances \cite{martinez2008multifunction}.
These delays may significantly complicate efforts to maintain coordination, control stability, and a seamless teleoperation experience \cite{lawrence1993stability}. They compromise the robot's duties' precision and safety as well as the natural feedback system between the remote robot and the human operator \cite{matheson2019human}. It is crucial to enable smooth exchange of contact force data from the slave robot to the master control system when controlling a robot remotely using a teleoperator. This communication serves the purpose of kinesthetically coupling the operator to the robot's environment, as mentioned in \cite{desoer2009feedback}, thereby significantly enhancing the operator's sense of telepresence and overall control \cite{anderson_bilateral}. However, if transmission delays are present, force feedback can have a destabilizing effect \cite{kim2004transatlantic}, \cite{fundamental_stability}.

Stability concerns are frequently resolved in some force-reflecting teleoperation systems by adding a significant amount of damping at various points across the system \cite{salcudean1995design}, \cite{kuchenbecker2004canceling}. However, this approach does not offer any formal stability guarantees and may considerably lower system performance \cite{main_ref}. The author in \cite{anderson_bilateral} discussed an approach to maintain stability in a force-reflecting bilateral teleoperator in the presence of a time delay. The stability of the transmission-delayed teleoperated system is examined in the works by \cite{lawrence1993stability}, \cite{reut_stability}, and \cite{hirai_stability}. Existing literature, see for example \cite{main_ref}, \cite{lawrence1993stability}, \cite{telerobotics}, \cite{Hokayem}, \cite{chan_teleop}, \cite{leung_teleop}, and \cite{anderson1992asymptotic}, addressed the stability of time-delayed teleoperation systems by using a passivity-based formalism. "Passivity-based formalism" is a systematic approach that employs mathematical models and control strategies to design teleoperation systems based on the concept of "passivity", which ensures system stability by regulating energy flow between the operator and the remote system. In the pursuit of a deeper comprehension of time-delayed transmissions and their intricate interplay with nonlinear dynamic systems, the author in \cite{main_ref} adeptly harnesses the concepts of wave variables and wave transmission to delve into system stability through the prism of passivity. This pioneering endeavor heralds a transformative design paradigm for teleoperation systems, ingeniously integrating communication pathways between two impedance controllers, thereby ensuring system stability as an outcome.

Nevertheless, the employment of wave variables and wave transmission within the system under higher time delay conditions exacerbates the convergence rate, inducing an undesirable oscillatory behavior, primarily attributed to the prolonged settling time of power variables convergence. This is demonstrated by the experimental results in the present work. \textbf{Contributions: }This paper presents a preliminary investigation into mitigating the impact of time delay and its destabilizing effect on force-reflecting telerobotic systems by using the concept of passivity-based formalism. Addressing this challenge, we propose a strategy within this paper to counteract the adverse effects of increased time delays on convergence. Specifically, we employ predictors such as the Smith predictor \cite{fundamental_smith, smith} and the minimum Jerk (MJ) predictor \cite{fundamental_mj, inbook} at the operator-end and remote-end, respectively, as part of our approach to enhance convergence performance. Through experimental validation, we vividly illustrated the disruptive influence of increased delays on convergence and underscored the effectiveness of mitigating this challenge through the integration of predictors. The experimentation is conducted on an F1/10 autonomous vehicle platform, teleoperated through a driving cockpit connected wirelessly to the platform and subjected to programmable network delays. 

The remainder of this paper is organized as follows: Section II contains notations employed throughout the paper. Section III contains preliminary information to briefly review the passivity formalism, wave variable and wave scattering concept, Smith predictor, and MJ predictor, which are required to follow the paper. Section IV presents the problem formulation with a theoretical analysis of the destabilizing effect of time delay. Section V presents the experimental and simulation results for the improved convergence of power variables for longer variable time delays, along with a comparison with existing results in the literature. The conclusion and future work are given in Section VI.

\section{Notation}

The following notions and conventions are employed throughout the paper:
$\mathbb{R},\mathbb{R}^n,\mathbb{R}^{n\times m}$  denote the space of real numbers, real vectors of length $n$ and real matrices of $n$ rows and $m$ columns respectively. 
$\mathbb{R}_+$ denotes positive real numbers. Normal-face lower-case letters ($x\in\mathbb{R}$) are used to represent real scalars, bold-face lower-case letter ($\mathbf{x}\in\mathbb{R}^n$) represents vectors, while normal-face upper case ($X\in\mathbb{R}^{n\times m}$) represents matrices. $\dot{x}$ is the time derivative of $x$. The set of all natural numbers is denoted by $\mathcal{N}$. $X^\top$ denotes the transpose of the quantity $X$. $\mathbf{x}^\top \mathbf{y}$ donates the scalar product of two vectors.

\section{Preliminaries}
\subsection{Passivity formalism}
This formalism employs energy principles to facilitate stability analysis, ensuring global stability in nonlinear systems \cite{desoer2009feedback}. Let $\psi$ represent the power (which may not correspond to actual physical power) entering a system, $\mathbf{x}$ be the input signal and $\mathbf{y}$ be the output signal of the system. Let $E$ be a lower-bounded energy function and $\zeta$ be a non-negative dissipation function. Then the system is said to be passive, as defined in \cite{slotine1991applied}, if it obeys the following condition

\begin{equation} \label{eq_1}
\psi = \mathbf{{x}}^\top \mathbf{y}= \frac{d E}{d t}+\mathbf{\zeta}
\end{equation}
This implies that the total energy supplied by the system up to time $t$ is limited to the initial stored energy:
\begin{align*}
\int_0^t \psi d \tau=\int_0^t (\mathbf{{x}}^\top \mathbf{y}) d \tau & =E(t)-E(0)+\int_0^t \zeta(\tau) d \tau
\end{align*}

\begin{figure}[t!]
\centerline{\includegraphics[scale=0.4]{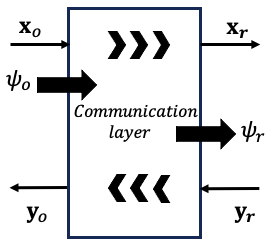}}
\caption{Passive 2-port element with operator-end and remote-end}
\label{fig1}
\end{figure}

If there is no power dissipation at any given point in time, the system is referred to as lossless. On the other hand, if power dissipation is positive as long as the stored energy has not reached its minimum threshold, then the system is considered strictly passive. Utilizing the stored energy as a Lyapunov-type function enables rapid analysis of stability and demonstrates that a passive system remains stable without any external input \cite{main_ref}. Asymptotic stability is achieved by a strictly passive system, provided that the stored energy is positively dependent on all states of the system. When two passive systems are connected in feedback or parallel configuration, then the overall system is also passive. This property is useful in force-reflecting teleoperation to manipulate an arbitrary passive environment without compromising global stability, which is a necessary condition for overall stability including the remote environment. Similarly, multiple passive 2-port elements can be cascaded into an overall passive 2-port element. 

For simplicity, we introduce a convention specific to 2-port elements, as depicted in Fig.\ref{fig1}, where power is designated to enter the system at the operator side as $\psi_o$ and leave the system at the remote side as $\psi_r$. Positive power is defined as entering the system. The total power flow is given below based on \eqref{eq_1}
$$ \psi = \psi_o - \psi_r = \mathbf{{x}}_o^\top \mathbf{y}_o - \mathbf{{x}}_r^\top \mathbf{y}_r $$

\subsection{Wave variable and Wave scattering}
The relationship between wave scattering and passivity is explored in detail in  \cite{main_ref}, while the practical implementation of these concepts can be located in  \cite{implement_wave}. At its core, this concept hinges on the notion that the power flow can be partitioned into the input and output parts associated with the input and output wave. Considering $\mathbf{u}_i$ and $\mathbf{v}_i$ as the input and output waves variables, respectively, where subscript $i$ is $o$ for the operator-end and $r$ for the remote-end, respectively. The total power flow as a function of input and output wave variables for a 2-port element, as depicted in Fig. \ref{fig1}, is
\begin{equation} \label{eq_2}
 \psi =  \frac{1}{2} \mathbf{u_o}^\top \mathbf{u_o}-\frac{1}{2} \mathbf{v_o}^\top \mathbf{v_o}+\frac{1}{2} \mathbf{u_r}^\top \mathbf{u_r}-\frac{1}{2} \mathbf{v_r}^\top \mathbf{v_r} 
\end{equation}

Based on \cite{main_ref}, the transformation between power variables $\mathbf{x}_i$ and $\mathbf{y}_i$ and wave variables $ \mathbf{u}_i $ and $ \mathbf{v}_i $ at operator-end and remote-end, respectively are described as
\begin{equation}
\begin{aligned} \label{eq_3} 
    \begin{bmatrix}
    \mathbf{u_o}\\\mathbf{v_o}
    \end{bmatrix} 
    = \begin{bmatrix}
        \beta & \alpha\\ -\beta & \alpha
        \end{bmatrix} 
        \begin{bmatrix}
    \mathbf{x_o}\\\mathbf{y_o}
    \end{bmatrix} , 
    \quad
    \begin{bmatrix}
    \mathbf{u_r}\\\mathbf{v_r}
    \end{bmatrix} 
    = \begin{bmatrix}
        -\beta & \alpha\\ \beta & \alpha
        \end{bmatrix} 
        \begin{bmatrix}
    \mathbf{x_r}\\\mathbf{y_r}
    \end{bmatrix}
\end{aligned}
\end{equation}
where $\alpha = (1/{\sqrt{2 b}})$ and $\beta = (\sqrt{b/2})$ and the strictly positive parameter $b$ is chosen based on the concept of wave impedance matching. The inverse transformation is give by: 
\begin{equation}
\begin{aligned} \label{eq_4} 
    \begin{bmatrix}
    \mathbf{x_o}\\\mathbf{y_o}
    \end{bmatrix} 
    = \begin{bmatrix}
        \alpha & -\alpha\\ \beta & \beta
        \end{bmatrix} 
        \begin{bmatrix}
    \mathbf{u_o}\\\mathbf{v_o}
    \end{bmatrix} , 
    \quad
    \begin{bmatrix}
    \mathbf{x_r}\\\mathbf{y_r}
    \end{bmatrix} 
    = \begin{bmatrix}
        -\alpha & \alpha\\ \beta & \beta
        \end{bmatrix} 
        \begin{bmatrix}
    \mathbf{u_r}\\\mathbf{v_r}
        \end{bmatrix}.
\end{aligned}
\end{equation}

\begin{figure}[t!]
\centerline{\includegraphics[scale=0.325]{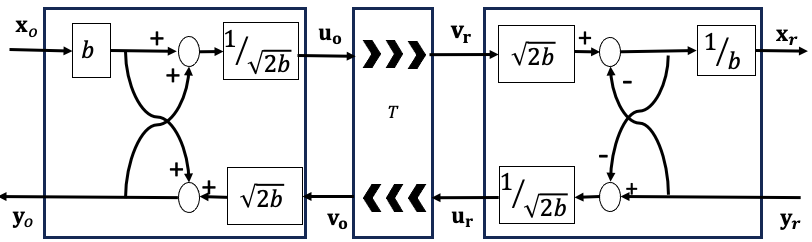}}
\caption{wave transformation between power and wave variables}
\label{fig3}
\end{figure}

The transformations are shown graphically in Fig. \ref{fig3} (assuming operator-side velocity $\mathbf{x}_o$ and remote-end force $\mathbf{y}_r$ are given).

The system is passive if the energy provided by output waves is limited to the energy received from input waves. Using \eqref{eq_2} the system is said to be passive if
$$
\int_0^t \frac{1}{2} (\mathbf{v}_o ^\top \mathbf{v}_o + \mathbf{v}_r ^\top \mathbf{v}_r) d \tau \leq \int_0^t \frac{1}{2} (\mathbf{u}_o ^\top \mathbf{u}_o + \mathbf{u}_r ^\top \mathbf{u}_r) d \tau.
$$

This is satisfied when the output wave amplitude is bounded by the amplitude of the delayed input wave even under time delay. We can therefore include arbitrary time delays into the system described by wave variables in a passive and hence stable fashion \cite{main_ref}.

\subsection{Smith-Predictor}
\begin{figure}[h!]
\centerline{\includegraphics[scale=0.35]{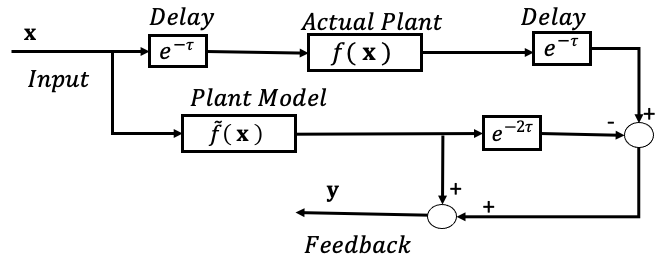}}
\caption{Schematic diagram of the Smith Predictor}
\label{fig4}
\end{figure}

The Smith Predictor \cite{fundamental_smith} addresses the challenge of instability and sub-optimal performance caused by feedback time delays by incorporating a predictive component within the control system. It hinges on the availability of a plant model, denoted as $\widetilde{f}(\mathbf{x})$, which characterizes the dynamic behavior of the system. The predictor anticipates feedback from the remote system by considering the given input while accounting for time delay $\tau$ \cite{smith}. The schematic diagram in Fig.\ref{fig4} illustrates the process, wherein the input $\mathbf{x}$ is provided to the predictor, yielding the predicted feedback $\mathbf{y}$ as the output. 
\begin{align*}
  \mathbf{y} & = f\left(\mathbf{x}\right)e^{-2 \tau}-\widetilde{f}(\mathbf{x})e^{-2 \tau}+\widetilde{f}(\mathbf{x}). 
\end{align*}
This predictive strategy aids in mitigating the effects of time delays, enhancing stability, and improving overall control system performance.

\subsection{MJ-Predictor}
Minimum Jerk Trajectory (MJT) \cite{inbook,fundamental_mj} between two points is the path that passes through the two points such that the integral of the jerk\footnote{jerk is the time derivative of acceleration} is minimized. The MJ-predictor is a prediction of future points based on the minimum trajectory interpolant of the historical data. Let $\mathbf{x}_1 = \mathbf{x}(t_1)$ and $\mathbf{x}_2=\mathbf{x}(t_2)$ be the observed points on the trajectorty $\mathbf{x}(t)$ at times $t_1$ and $t_2$ respectively. Then, the MJT between $\mathbf{x}_1$ and $\mathbf{x}_2$ is given by:
\begin{align}\label{eq_mj}
    \hat{\mathbf{x}}(t) = \mathbf{x}_1 +\left(\mathbf{x}_2\ -\mathbf{x}_1\right)\ (6\gamma^5\ -\ 15\gamma^4\ +\ 10\gamma^3),
\end{align}
where $\gamma = (t-t_1)/(t_2-t_1)$. Based on this, future values of $\mathbf{x}(t)$ are predicted as $\hat{\mathbf{x}}(t)$ for $t>t_2$.

\iffalse
The MJ-Predictor is used to predict the minimum jerk of any movement using the concept of the  Minimum Jerk Trajectory (MJT). The MJT is a trajectory designed to minimize the rate of change of acceleration (jerk), which is adopted by the author of \cite{inbook} from \cite{fundamental_mj} to derive the following equation \eqref{eq_mj}, which is used to calculate the jerk trajectory for the interval $t_1$ to $t_2$, given two inputs $\theta_1$ and $\theta_2$ at times $t_1$ and $t_2$, respectively. 
% There are other MJ predictors \cite{prev_work1_mj} and \cite{prev_work2_mj}, which utilize the concept for MJT prediction
\begin{equation} \label{eq_mj}
    \theta(\gamma)\ =\ \theta_1\ +\ (\theta_2\ -\ \theta_1)\ (6\gamma^5\ -\ 15\gamma^4\ +\ 10\gamma^3)
\end{equation} 

where $\gamma = (t-t_1)/(t_2-t_1)$ and $\theta(\gamma)$ is the trajectory point at time $t$. 

Let $\theta_0 = 0$ at $t_1=0$ and $\theta(\gamma) = \theta_p$, which is a known value at time $t=t_p$ then the above equation \eqref{eq_mj} can be rewritten as 
\begin{equation} \label{eq_mj2}
\theta_2 =\theta_p / (6\gamma^5\ -\ 15\gamma^4\ +\ 10\gamma^3)
\end{equation}
Where $\gamma = t_p/t_f $ and $t_f$ is the time at which the minimum jerk trajectory point is predicted.
\fi

\section{Problem Formulation}
The communications element in teleoperation connects local and remote systems to close the control loop by transmitting data to and from both sites. It introduces time delays, which may be caused by physical transmission times or communication bandwidth limitations. Time delays can cause instability in feedback systems, including in force-reflecting teleoperation due to communication delays between local and remote sites. If the given operator and remote sites are passive and if the communication element is passive then the entire system is stable. Using the fact that the connected multiple passive subsystems result in a combined passive system, we can separate the communication element, as shown in Fig.\ref{fig1} (described by a 2-port element, with one connected to the operator site and the other to the remote system), to study the instability caused by the time delay. Considering $\tau$ as the time delay in the communication layer, we have
\begin{equation} \label{eq_5}
    \mathbf{x}_r(t)= \mathbf{x}_o(t-\tau), \quad \mathbf{y}_o(t) = \mathbf{y}_r(t-\tau)
\end{equation} 
Now, substituting \eqref{eq_3} into \eqref{eq_2} to look at the power flow into the system, we get the following: 
\begin{equation} \label{eq_6}
\begin{aligned}
\psi =  & \frac{1}{2 b} \mathbf{y}_o^2(t)+\frac{b}{2} \mathbf{x}_o^2(t)-\frac{1}{2 b}\left(\mathbf{y}_o(t)-b \mathbf{x}_o(t)\right)^2 \\
& +\frac{1}{2 b} \mathbf{y}_r^2(t)+\frac{b}{2} \mathbf{x}_r^2(t)-\frac{1}{2 b}\left(\mathbf{y}_r(t)+b \mathbf{x}_r(t)\right)^2.
\end{aligned}
\end{equation}
Using \eqref{eq_5}, \eqref{eq_6} can be simplified to 
\begin{equation} \label{eq_7}
\begin{aligned}
\psi & = \frac{1}{b} \mathbf{y}_o^2(t) - \frac{1}{2 b}\left(\mathbf{y}_o(t)-b \mathbf{x}_o(t)\right)^2 -\frac{1}{2 b}\left(\mathbf{y}_r+b \mathbf{x}_r\right)^2(t) \\
& + b \mathbf{x}_r^2(t) + \frac{d}{d t} \int_{t-\tau}^t \left( \frac{b}{2} \mathbf{x}_o^2(\tau)+\frac{1}{2 b} \mathbf{y}_r^2(\tau) \right) d \tau.
\end{aligned}
\end{equation}
Comparing \eqref{eq_1} with \eqref{eq_7}, we get the energy storage function, $E$ and power dissipation, $\zeta$ as
\begin{align*}
E & = \int_{t-\tau}^t \left( \frac{b}{2} \mathbf{x}_o^2(\tau)+\frac{1}{2 b} \mathbf{y}_r^2(\tau) \right) d \tau
\end{align*}
and
\begin{align*}
\zeta = & \frac{1}{b} \mathbf{y}_o^2(t) - \frac{1}{2 b}\left(\mathbf{y}_o(t)-b \mathbf{x}_o(t)\right)^2 \\
& + b \mathbf{x}_r^2(t) -\frac{1}{2 b}\left(\mathbf{y}_r(t)+b \mathbf{x}_r(t)\right)^2. 
\end{align*}
For the communication element to be passive, the power dissipation $\zeta$ must be positive. However for certain input values of power variables $\mathbf{x}_o$ and $\mathbf{y}_r$, the power dissipation, $\zeta$, can be negative, which can be seen from the above expression. It is highly unfavorable to employ non-passive communications, which leads to instability. Therefore, the communications layer should be passive. 

In standard communication setups, energy injection can cause instability, often resolved by adding ample damping. The increased dissipation ensures energy absorption, rendering the modified communication system stable. However, this method changes the dynamics and introduces unwanted effects and continuous power input is required to sustain a constant motion and constant force reflection. This approach may not be suitable for direct applications.

A better solution is to transmit the wave variables instead of the power variables. In the previous section, we showed that the system can be unstable when using power variables $\mathbf{x}_0$ and $\mathbf{y}_0$, instead using wave variables makes the communication element passive, which is shown below. Considering the setup shown in Fig.\ref{fig3} with a time delay $\tau$ in the communication layer, we have
$$
\mathbf{v_o}(t)=\mathbf{u_r}(t-\tau), \quad \mathbf{v_r}(t)=\mathbf{u_o}(t-\tau).
$$
Then, the power flow into the setup is given by
\begin{align*}
\psi & =\frac{1}{2} \mathbf{u_o}(t)^2-\frac{1}{2} \mathbf{v_o}(t)^2+\frac{1}{2} \mathbf{u_r}(t)^2-\frac{1}{2} \mathbf{v_r}(t)^2 \\
& =\frac{1}{2} \mathbf{u_o}(t)^2-\frac{1}{2} \mathbf{u_r}(t-\tau)^2+\frac{1}{2} \mathbf{u_r}(t)^2-\frac{1}{2} \mathbf{u_o}(t-\tau)^2 \\
& =\frac{d}{d t}\left[\int_{t-\tau}^t \frac{1}{2} \mathbf{u_o}(\tau)^2+\frac{1}{2} \mathbf{u_r}(\tau)^2 d \tau\right].
\end{align*}
Here, we can see that the lossless passive communication is achieved using wave variables and overall stability is preserved. However, using a wave transmission scheme can also introduce wave reflections, which can be avoided by ensuring the impedance of the wave transmission $b$ is aligned with the rest of the system, either by choice of parameter or by including additional termination elements. Using the additional termination elements to the setup shown in Fig.\ref{fig3}, the modified setup will be the setup shown in Fig.\ref{fig4}, using which the wave transformations are then governed by 
\begin{equation}
\begin{aligned} \label{eq_8}
    & \mathbf{u_o} =\sqrt{\frac{b}{2}} \mathbf{x}_o,  \quad \mathbf{u_r}  =\frac{1}{\sqrt{2 b}} \mathbf{y}_r \\
    & \mathbf{y}_o=\frac{b}{2} \mathbf{x}_o+\sqrt{\frac{b}{2}} \mathbf{v_o}, \quad \mathbf{x}_r=-\frac{1}{2 b} \mathbf{y}_r+\frac{1}{\sqrt{2 b}} \mathbf{v_r}.
\end{aligned}
\end{equation}
\begin{figure}[t!]
\centerline{\includegraphics[scale=0.27]{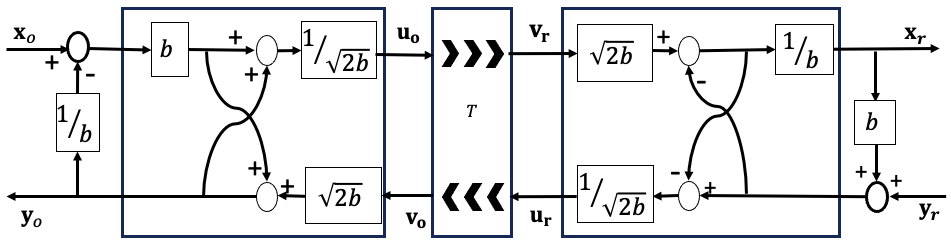}}
\caption{wave transformation between power and wave variables after impedance matching}
\label{fig5}
\end{figure}
This setup solves the issue of instability by making the overall system passive and avoiding wave reflections. However, for higher time delays the power variables do not converge and have an oscillatory behavior, which is shown in the experimental results in the experimental section. To avoid this behavior wave variables have to be predicted like \cite{wave_predictor}, in this paper, two different predictors are used: the Smith predictor to predict the remote wave variables and the MJ predictor to predict the operator-end wave variables.

\begin{figure}[b!]
\centerline{\includegraphics[scale=0.3]{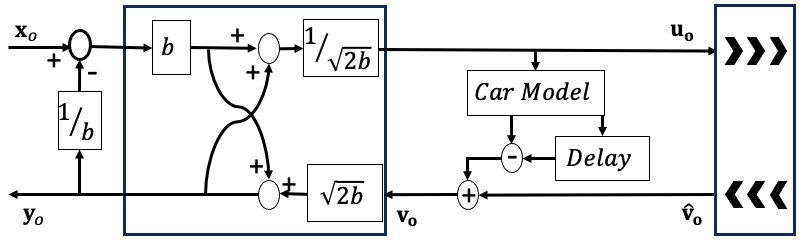}}
\caption{Wave Variable Prediction - Operator End with Smith Predictor}
\label{fig6}
\end{figure}

The above-mentioned setup along with the predictors is depicted in Fig.\ref{fig6} from which it can be observed that the Smith predictor takes the wave variable $\mathbf{u_o}(t)$ as the input and gives the output wave variable from the plant model at time $t$ and $t-2\tau$. Let the plant model gives wave variable output be $\mathbf{{y}}(t)$ and $\mathbf{{y}}(t-2\tau)$ at time $t$ and $t-2\tau$, respectively, then the wave variable $\mathbf{v_o}(t)$ can be expressed as follows
\begin{equation}  \label{eq_9}
\begin{aligned}
        \mathbf{v_o}(t) & = \mathbf{u_r}(t-\tau) + \mathbf{{y}}(t) - \mathbf{{y}}(t-2\tau) \\
    & = \mathbf{\hat{v}_o}(t) + \mathbf{{y}}(t) - \mathbf{{y}}(t-2\tau),
\end{aligned}
\end{equation}
where $\mathbf{u_r}(t-\tau)$ is delayed wave variable feedback from the remote side, which is a known value. Considering this fact, equation \eqref{eq_9} and the expression from \eqref{eq_8}, we get the feedback at the operator end at time $t$, $\mathbf{y}_o(t)$ as follows
\begin{equation}  \label{eq_10}
     \mathbf{y}_o(t) = \frac{b}{2} \mathbf{x}_o(t)+\sqrt{\frac{b}{2}} [\mathbf{\hat{v}_o}(t) + \mathbf{{y}}(t) - \mathbf{{y}}(t-2\tau)]. 
\end{equation}
The feedback at the operator end at $t$ can be found using the above expression since the $\mathbf{x}_o(t)$ and $\hat{\mathbf{v_o}}$ are the known values.

\begin{figure}[t!]
\centerline{\includegraphics[scale=0.34]{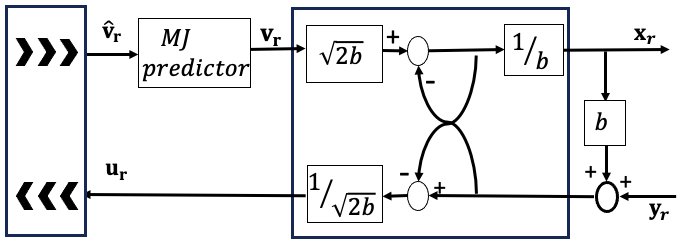}}
\caption{Wave Variable Prediction - Remote End with MJ predictor}
\label{fig7}
\end{figure}

Smith predictor is used on the operator end to predict the remote end wave variable since the model of the remote end is considered available but the model of the operator end is considered unknown so we use the MJ predictor. At the time of delay remote end doesn't receive any input wave variable so to predict human command MJ predictor is used. It predicts the wave variable $\mathbf{v_r}$ from the operator end by predicting the trajectory of the $\mathbf{\hat{v}_o}$ wave variable, as shown in Fig.\ref{fig6} \& Fig.\ref{fig7}. This prediction provides the data to have continued communication even if the delay is high by giving the predicted values instead of using the same previous values. 

The rationale for employing two distinct simulators lies in our intention to forecast the behaviors of two different systems: one system being operated by humans and the other by an autonomous entity. Given the availability of comprehensive kinematic and dynamic models for the autonomous system, the utilization of the Smith predictor is deemed advantageous, offering heightened precision in prediction. Conversely, for forecasting human motion over subsequent intervals, the MJ predictor emerges as the preferred choice due to the absence of a perfect human model. It is worth noting that employing the MJ predictor for the former system is suboptimal, as it primarily forecasts jerk, whereas the Smith predictor possesses the capability to anticipate the vehicle's motion more comprehensively. By leveraging these distinct predictors judiciously, we aim to enhance the predictive accuracy and effectiveness of our system across diverse operational scenarios.

\section{Implementation and Results}
The operator controls the remote-controlled car using a G29 cockpit controller. Feedback from the car is conveyed to the operator in two ways: the car's velocity and force (torque reflected on steering) are obtained from the sensor data provided by the vehicle. The communication infrastructure bridging the car and the cockpit is implemented by leveraging the Robot Operating System (ROS) framework.

\subsection{Simulation results}

\begin{figure}[h!]
\centerline{\includegraphics[scale=0.45]{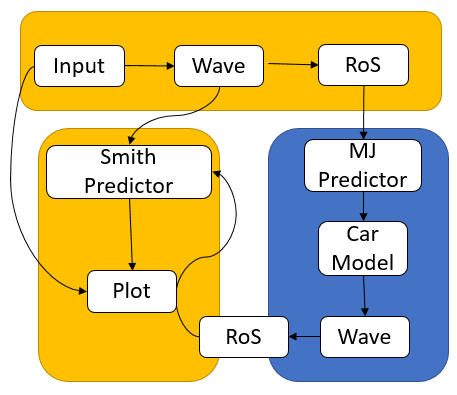}}
\caption{Communication architecture for the simulation environment.}
\label{fig12}
\end{figure}

For the simulation model, we used Simulink with ROS framework. The Simulink environment incorporates the ROS layer to faithfully replicate the configuration presented in Fig.\ref{fig12}. In the simulation, a step signal with an amplitude of 0.5 is considered as the input signal, $\mathbf{x_o}$, given to the model shown in Fig.\ref{fig12} for a duration of 10 seconds. Two distinct cases are considered during the simulation: one without any forced delay and another with a forced delay of 1 second. 

\begin{figure}[h!]
\centerline{\includegraphics[scale=0.45]{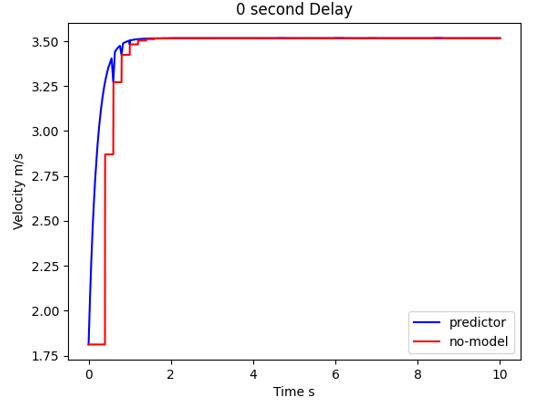}}
\caption{Simulation response of with- and without-predictors: 0sec delay}
\label{fig13}
\end{figure}

\begin{figure}[h!]
\centerline{\includegraphics[scale=0.45]{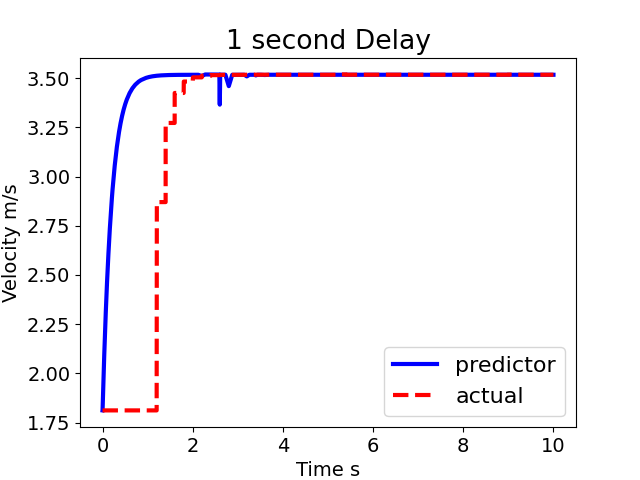}}
\caption{{Simulation response of with- and without-predictors: 1sec delay}}
\label{fig14}
\end{figure}

For both cases, the data is collected to examine and visualize and then plotted as shown in Fig. \ref{fig13} and Fig. \ref{fig14} from the observation of these figures we can find that even in the absence of forced delay, the actual system exhibits a temporal response, signifying that it requires a certain duration to converge. Conversely, under the influence of forced delay, the settling time significantly extends. \par

Importantly, for the model with the predictor, the settling time of the velocity response is less compared to the model without the predictor and the performance of the velocity response remains consistent for both the cases, demonstrating a stable response without any oscillatory behavior. For the predictor's case results, the consistent response showcases the predictor's effectiveness in mitigating delays and ensuring a prompt and reliable system response, irrespective of the temporal conditions imposed. For the simulation, the analysis is performed for the linear velocity but the same can be used for the angular velocity analysis as well.\par

% \begin{figure}[htbp]
% \centerline{\includegraphics[scale=0.35]{plot_simu.jpg}}
% \caption{Simulation Result}
% \label{fig14}
% \end{figure}

\subsection{Experimental validation}

\begin{figure}[h!]
\centerline{\includegraphics[scale=0.35]{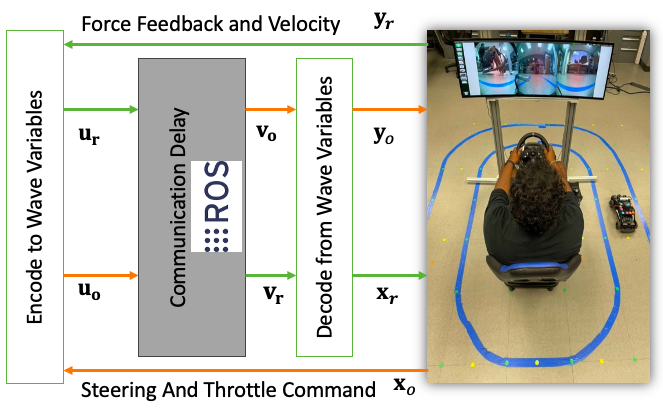}}
\caption{Experimental Setup with ROS communication, driving cockpit with video and steering torque feedback, and remote F1/10th autonomous vehicle.}
\label{fig8}
\end{figure}

The operator supplies throttle and steering angle as inputs $\mathbf{x_o}$ and receives feedback in the form of linear and angular velocities $\mathbf{y_o}$ as shown in Fig. \ref{fig8}. As depicted in Fig. \ref{fig12}, the operator's input undergoes encoding into wave variables before transmission over the Robot Operating System (ROS) network.

In the context of our simulation, we acknowledge the inclusion of a constant time delay. However, in our experimental setup, we aim to ascertain and quantify the teleoperating delay at both ends. Subsequently, we plan to dynamically reconfigure the predictor programmatically to align with the observed time delay. 

For our experimental setup, with an input range $\mathbf{x_o}$ from [-1, 1], the corresponding feedback $\mathbf{y_o}$ ranges from [7.5, 7.5], and in the absence of system delay, when $\mathbf{x_0}(t) = 1$, then $\mathbf{y_0}(t) = \mathbf{y_r}(t)$ = 7.5. By using \eqref{eq_8} and \eqref{eq_4} we get the impedance $b$ value as 7.5.

However, taking into account the influence of time delay contemplate $b$ to 8. This encoded input is forwarded to the MJ predictor, wherein the minimum jerk of the human operator's input is predicted, and subsequently, the wave variable is decoded. In instances of delay, the controller will not receive any input so it uses the previously received value and forecasts the current value, which is then transmitted to the car. The feedback from the car is encoded into a wave variable and seamlessly transmitted over ROS, with subsequent decoding transpiring at the operator's end, this process is depicted in Fig. \ref{fig8}.

The Smith predictor initiates its operations upon the commencement of operator input, regardless of immediate feedback from the car. Real-time calculation of communication delays is consistently updated to the predictor. Consequently, the Smith predictor promptly commences forecasting the car's response based on the provided input. Employing a technical approach, the predictor predicts the response not only for the current input but also factors in the time-delayed input from the history. This dual-input prediction mechanism enables the Smith predictor to finely predict the feedback, enhancing its effectiveness in compensating for system delays. This process has been carried out and the data is collected for three scenarios – without wave variables and predictors, with wave variables and without the predictor, and with wave variables and predictors, which are considered for three different cases.

\subsubsection{Case 1: Without forced delay}
\begin{figure}[h!]
\centerline{\includegraphics[scale=0.45]{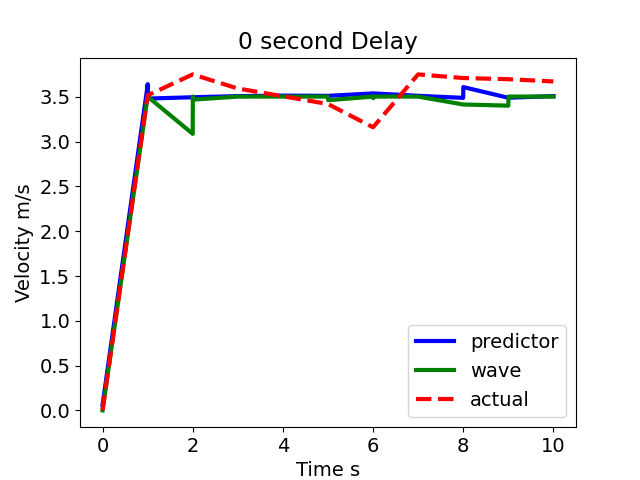}}
\caption{Experimental result with and without  predictors: No transmission delay}
\label{fig9}
\end{figure}

In the absence of any imposed delays within the system, indicating a 0-delay communication layer, a step input with an amplitude equal to half of the maximum allowable input is applied. Mathematically, the equilibrium point for this specified input is determined to be 3.5 m/s. Verification of this analytical result is depicted in Fig. \ref{fig9}. Notably, despite the absence of forced delays, the inherent natural delay of the system causes signal oscillation. However, this oscillation is effectively mitigated in a wave variable plot, showcasing convergence. Furthermore, predictors demonstrate commendable convergence, outperforming the wave variable controller in achieving stable responses.

\subsubsection{Case 2: With 500ms of forced delay}
\begin{figure}[h!]
\centerline{\includegraphics[scale=0.45]{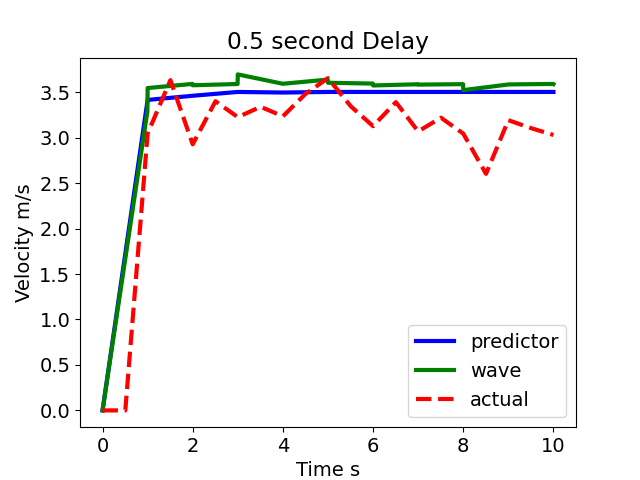}}
\caption{Experimental result with and without  predictors: 0.5s transmission delay.}
\label{fig10}
\end{figure}

We proceed to cross-examine the performance of the predictor and the wave variable controller under the influence of time delay, as illustrated in Fig. \ref{fig10}. The actual system exhibits oscillations around the equilibrium point without converging to it. In contrast, both the predictor and the wave variable controller demonstrate convergence to the desired point, albeit with a slightly longer settling time for the wave variable, which does not pose a critical issue; however, further implications arise when the transmission delay increases, and these will be explored in the next case.

\begin{figure}[h!]
\centerline{\includegraphics[scale=0.45]{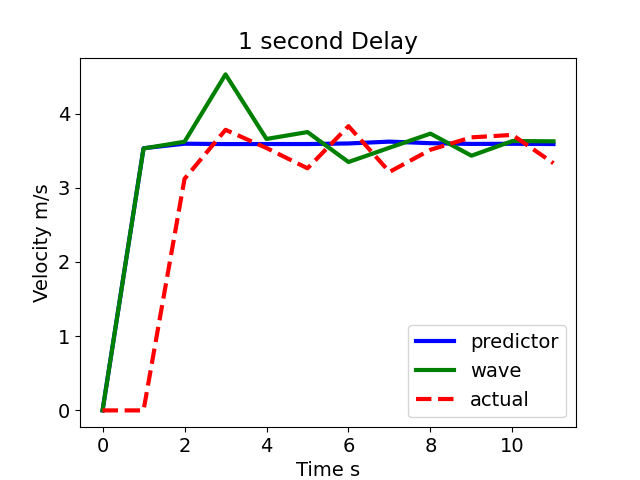}}
\caption{Experimental result with and without  predictors: 1s transmission delay}
\label{fig11}
\end{figure}

\subsubsection{Case 3: With 1s of forced delay}
In Fig. \ref{fig11}, observations indicate that for both the scenarios with and without the wave variables, the velocity response exhibits oscillations around the equilibrium point. Whereas, the predictor controller converges with the same settling time as observed in case 2, which is evidenced by the results. This implies that, in the presence of increased delay, the predictor controller maintains its ability to converge efficiently without increasing the settling time, showcasing its ability to outperform the default system and the system with the wave variable controller.

For the experiment, the analysis is performed for the linear velocity but the same can be used for the angular velocity analysis as well by commanding the steering angle and velocity of the car/robot.

\section{Conclusion }
In conclusion, this research paper has showcased a successful integration of both the Smith predictor and MJ predictor into a system, resulting in significant enhancements to the convergence rate and a substantial reduction in the settling time required to reach a steady state, even under conditions of higher transmission delays. Through experimental validation, we have vividly demonstrated how the sole reliance on wave variables can impede the convergence of power variables for higher delays, whereas the use of predictors effectively resolves this issue. By introducing these predictors, we have not only ensured the overall stability of the system but also improved the convergence rate, regardless of the extent of time delays. This achievement marks a significant milestone in facilitating reliable communication between human operators and remote-controlled systems, thus paving the way for more robust and efficient teleoperation across a wide range of applications.

\bibliographystyle{IEEEtran}
\bibliography{References}

\end{document}